\documentclass[11pt]{amsart}
\usepackage{latexsym,graphicx}

\numberwithin{equation}{section}
\theoremstyle{plain}

\theoremstyle{remark}

\theoremstyle{definition}

\newcommand{\D}{{\mathcal D}}
\newcommand{\E}{\mathcal E}

\newcommand{\G}{{\mathcal G}}

\newcommand{\K}{{\mathcal K}}
\renewcommand{\L}{{\mathcal L}}
\newcommand{\M}{{\mathcal M}}
\newcommand{\N}{\mathbb N}

\newcommand{\V}{{\mathcal V}}

\newcommand{\dist}{\operatorname{dist}}

\newcommand{\fp}{\operatorname{FP}}

\newcommand{\Int}{\operatorname{Int}}

\renewcommand{\span}{\operatorname{span}}

\newcommand{\supp}{\operatorname{Supp}}

\def\half{{1 \over 2}}
\def\la{\lambda}

\newcommand{\oa}{\overrightarrow}

\newcommand{\ol}{\overline}

\def\XXint#1#2#3{{\setbox0=\hbox{$#1{#2#3}{\int}$}
      \vcenter{\hbox{$#2#3$}}\kern-.5\wd0}}

\address{ Dept. of Mathematics, Rice  University, 6100 Main St., Houston, 77005 TX, U.S.A.
\\ {\sl E-mail address:}  {\bf harvey@rice.edu}}

\address{ Department of Mathematics, Stony Brook University, Stony Brook, NY 11790, U.S.A.
\\ {\sl E-mail address:}  {\bf blaine@math.sunysb.edu}}

\address{  Institute of Mathematics, Cracow University of Technology, Warszawska 24, 31-155
    Krak\'{o}w, Poland
\\ {\sl E-mail address:}  {\bf splis@pk.edu.pl}}

\begin{document}

\def\cal{\mathcal}

\font\tpt=cmr10 at 12 pt
\font\fpt=cmr10 at 14 pt

\font \fr = eufm10



\overfullrule=0in

\def\boxit#1{\hbox{\vrule
 \vtop{%
  \vbox{\hrule\kern 2pt %
     \hbox{\kern 2pt #1\kern 2pt}}%
   \kern 2pt \hrule }%
  \vrule}}

  \def\harr#1#2{\ \smash{\mathop{\hbox to .3in{\rightarrowfill}}\limits^{\scriptstyle#1}_{\scriptstyle#2}}\ }

\def\ALEX{1}
\def\AGV{2}
\def\ASSA{3}
\def\ASSB{4}
\def\BTA{5}
\def\BTB{6}
\def\BTC{7}
\def\CLN{8}
\def\CIL{9}
\def\CRA{10}
\def\GAR{11}
\def\HAR{12}
\def\DDD{13}
\def\PUP{14}
\def\DDR{15}
\def\HYP{16}
\def\HLGAR{17}
\def\REST{18}
\def\AC{19}
\def\pPSH{20}
\def\SURVEY{21}
\def\BELL{22}
\def\ASPECTS{23}
\def\HP{24}
\def\HPP{25}
\def\LAB{26}
\def\LLAB{27}
\def\LLLAB{28}
\def\LAN{29}
\def\LEL{30}
\def\POGA{31}
\def\POGB{32}
\def\SHF{33}
\def\TWA{34}
\def\TWB{35}
\def\TWC{36}
\def\VIT{37}

 \def\GG{{{\bf G} \!\!\!\! {\rm l}}\ }

\def\GL{{\rm GL}}

\def\bll{I \!\! L}

\def\bra#1#2{\langle #1, #2\rangle}
\def\bbf{{\bf F}}
\def\bbj{{\bf J}}
\def\Jtn{{\bbj}^2_n}  \def\JtN{{\bbj}^2_N}  \def\JoN{{\bbj}^1_N}
\def\jt{j^2}
\def\jtx{\jt_x}
\def\Jt{J^2}
\def\Jtx{\Jt_x}
\def\bpp{{\bf P}^+}
\def\bpt{{\wt{\bf P}}}
\def\fsh{$F$-subharmonic }
\def\mo{monotonicity }
\def\jet{(r,p,A)}
\def\ss{\subset}
\def\sse{\subseteq}
\def\half{\hbox{${1\over 2}$}}
\def\smfrac#1#2{\hbox{${#1\over #2}$}}
\def\oa#1{\overrightarrow #1}
\def\dim{{\rm dim}}
\def\dist{{\rm dist}}
\def\codim{{\rm codim}}
\def\deg{{\rm deg}}
\def\rank{{\rm rank}}
\def\log{{\rm log}}
\def\Hess{{\rm Hess}}
\def\Hessyp{{\rm Hess}_{\rm SYP}}
\def\trace{{\rm trace}}
\def\tr{{\rm tr}}
\def\max{{\rm max}}
\def\min{{\rm min}}
\def\span{{\rm span\,}}
\def\Hom{{\rm Hom\,}}
\def\det{{\rm det}}
\def\End{{\rm End}}
\def\Sym{{\rm Sym}^2}
\def\diag{{\rm diag}}
\def\pt{{\rm pt}}
\def\Spec{{\rm Spec}}
\def\pr{{\rm pr}}
\def\Id{{\rm Id}}
\def\Grass{{\rm Grass}}
\def\Herm#1{{\rm Herm}_{#1}(V)}
\def\arr{\longrightarrow}
\def\supp{{\rm supp}}
\def\Link{{\rm Link}}
\def\Wind{{\rm Wind}}
\def\Div{{\rm Div}}
\def\vol{{\rm vol}}
\def\foral{\qquad {\rm for\ all\ \ }}
\def\fpsh{{\cal PSH}(X,\f)}
\def\Core{{\rm Core}}
\def\dis{f_M}
\def\Re{{\rm Re}}
\def\rn{\bbr^n}
\def\pp{\cp^+}
\def\plp{\cp_+}
\def\Int{{\rm Int}}
\def\cix{C^{\infty}(X)}
\def\Gr#1{G(#1,\rn)}
\def\Symn{{\Sym(\rn)}}
\def\SymN{{\Sym(\bbr^N)}}
\def\Gpn{G(p,\rn)}
\def\fd{{\rm free-dim}}
\def\SA{{\rm SA}}
 \def\cd{{\cal C}}
 \def\cdt{{\widetilde \cd}}
 \def\cm{{\cal M}}
 \def\cmt{{\widetilde \cm}}

\def\Theorem#1{\medskip\noindent {\bf THEOREM \bf #1.}}
\def\Prop#1{\medskip\noindent {\bf Proposition #1.}}
\def\Cor#1{\medskip\noindent {\bf Corollary #1.}}
\def\Lemma#1{\medskip\noindent {\bf Lemma #1.}}
\def\Remark#1{\medskip\noindent {\bf Remark #1.}}
\def\Note#1{\medskip\noindent {\bf Note #1.}}
\def\Def#1{\medskip\noindent {\bf Definition #1.}}
\def\Claim#1{\medskip\noindent {\bf Claim #1.}}
\def\Conj#1{\medskip\noindent {\bf Conjecture \bf    #1.}}
\def\Ex#1{\medskip\noindent {\bf Example \bf    #1.}}
\def\Qu#1{\medskip\noindent {\bf Question \bf    #1.}}
\def\Exercise#1{\medskip\noindent {\bf Exercise \bf    #1.}}

\def\HoQu#1{ {\AAA T\BBB HE\ \AAA H\BBB ODGE\ \AAA Q\BBB UESTION \bf    #1.}}

\def\pf{\medskip\noindent {\bf Proof.}\ }
\def\qed{\hfill  $\vrule width5pt height5pt depth0pt$}
\def\equdef{\buildrel {\rm def} \over  =}
\def\qedqed{\hfill  $\vrule width5pt height5pt depth0pt$ $\vrule width5pt height5pt depth0pt$}
\def\mathqed{  \vrule width5pt height5pt depth0pt}

\def\V{W}

\def\df{d^{\phi}}
\def\hk{\_{\rm l}\,}
\def\n{\nabla}
\def\w{\wedge}

\def\cu{{\cal U}}   \def\cc{{\cal C}}   \def\cb{{\cal B}}  \def\cz{{\cal Z}}
\def\cv{{\cal V}}   \def\cp{{\cal P}}   \def\ca{{\cal A}}
\def\cw{{\cal W}}   \def\co{{\cal O}}
\def\ce{{\cal E}}   \def\ck{{\cal K}}
\def\ch{{\cal H}}   \def\cm{{\cal M}}
\def\cs{{\cal S}}   \def\cn{{\cal N}}
\def\cd{{\cal D}}
\def\cl{{\cal L}}
\def\cp{{\cal P}}
\def\cf{{\cal F}}
\def\ccr{{\cal  R}}

\def\gerG{{\fr{\hbox{g}}}}
\def\gerB{{\fr{\hbox{B}}}}
\def\gerR{{\fr{\hbox{R}}}}
\def\p#1{{\bf P}^{#1}}
\def\vf{\varphi}

\def\wt{\widetilde}
\def\wh{\widehat}

\def\and{\qquad {\rm and} \qquad}
\def\arr{\longrightarrow}
\def\ol{\overline}
\def\bbr{{\mathbb R}}\def\bbh{{\mathbb H}}\def\bbo{{\mathbb O}}
\def\bbc{{\mathbb C}}
\def\bbq{{\mathbb Q}}
\def\bbz{{\mathbb Z}}
\def\bbp{{\mathbb P}}
\def\bbd{{\mathbb D}}

\def\a{\alpha}
\def\b{\beta}
\def\d{\delta}
\def\e{\epsilon}
\def\f{\phi}
\def\g{\gamma}
\def\k{\kappa}
\def\la{\lambda}
\def\o{\omega}

\def\s{\sigma}
\def\x{\xi}
\def\z{\zeta}

\def\D{\Delta}
\def\L{\Lambda}
\def\G{\Gamma}
\def\O{\Omega}

\def\bd{\partial}
\def\bdf{\partial_{\f}}
\def\lag{Lagrangian}
\def\psh{plurisubharmonic }
\def\ph{pluriharmonic }
\def\pph{partially pluriharmonic }
\def\omp{$\omega$-plurisubharmonic \ }
\def\ffl{$\f$-flat}
\def\PH#1{\widehat {#1}}
\def\lloc{L^1_{\rm loc}}
\def\dbar{\ol{\partial}}
\def\lp{\Lambda_+(\f)}
\def\lpp{\Lambda^+(\f)}
\def\bo{\partial \Omega}
\def\Ob{\overline{\O}}
\def\fc{$\phi$-convex }
\def\PSH{{ \rm PSH}}
\def\SH{{\rm SH}}
\def\totr{ $\phi$-free }
\def\BM{\lambda}
\def\Der{D}
\def\CH{{\cal H}}
\def\RH{\overline{\ch}^\f }
\def\pconv{$p$-convex}
\def\MA{MA}
\def\lagpsh{Lagrangian plurisubharmonic}
\def\hermsk{{\rm Herm}_{\rm skew}}
\def\PSHl{\PSH_{\rm Lag}}
 \def\ppsh{$\pp$-plurisubharmonic}
\def\fp{$\pp$-plurisubharmonic }
\def\fh{$\pp$-pluriharmonic }
\def\Symn{\Sym(\rn)}
 \def\ci{C^{\infty}}
\def\USC{{\rm USC}}
\def\fa{{\rm\ \  for\ all\ }}
\def\ppc{$\pp$-convex}
\def\cpt{\wt{\cp}}
\def\ft{\wt F}
\def\ob{\overline{\O}}
\def\Be{B_\e}
\def\K{{\rm K}}

\def\M{{\bf M}}
\def\N#1{C_{#1}}
\def\ds{Dirichlet set }
\def\dir{Dirichlet }
\def\Fa{{\oa F}}
\def\TR{{\cal T}}
 \def\LAG{{\rm LAG}}
 \def\ISO{{\rm ISO_p}}
 \def\Span{{\rm Span}}

\def\AA{1}
\def\BB{2}
\def\CC{3}
\def\DD{4}
\def\EE{5}
\def\FF{6}
\def\GGG{7}
\def\HH{8}
\def\II{9}
\def\JJ{10}
\def\KK{11}
\def\LL{12}
\def\MM{13}

\vskip .4in

\def\BLO{1}
\def\CRA{2}
\def\CIL{3}
\def\DS{4}
\def\DK{5}
\def\Fo{6}
\def\FW{7}
\def\DDD{8}
\def\PUP{9}
\def\DDR{10}
\def\GPSH{11}
\def\AC{12}
\def\Survey{13}
\def\Lu{14}
\def\LN{15}
\def\Pa{16}
\def\Plis{17}
\def\Pl{18}
\def\Pll{19}
\def\Plll{20}

\def\PSF{{\mathcal PSH}}

\def\E{E}
\def\bL{{\bf \Lambda}}
\font\headfont=cmr10 at 14 pt

\vskip .1in


\title[SMOOTH APPROXIMATION 
ON ALMOST COMPLEX MANIFOLDS]{SMOOTH APPROXIMATION OF \\PLURISUBHARMONIC FUNCTIONS\\
ON ALMOST COMPLEX MANIFOLDS
}
 
\date{\today}
\author{ F. Reese Harvey,  H. Blaine Lawson, Jr.
and Szymon Pli\'s}
 \thanks{The second author was partially supported by the NSF and IHES,
 and the third author
was partially supported by the NCN grants 2011/01/D/ST1/04192 and 2013/08/A/ST1/00312.
}

\maketitle

\centerline{\bf Abstract}
  \font\abstractfont=cmr10 at 10 pt
  
  {{\parindent= .2in\narrower \noindent
  
This note establishes smooth approximation from above for $J$-plurisubharmonic functions
on an almost complex manifold $(X,J)$.  The following theorem is proved.
Suppose $X$ is $J$-pseudoconvex, i.e., $X$ admits a smooth strictly $J$-plurisubharmonic exhaustion function.
Let $u$ be an (upper semi-continuous)  $J$-plurisubharmonic  function
on $X$. Then there exists a  sequence $u_j \in  C^\infty(X)$ of smooth strictly 
$J$-plurisubharmonic functions point-wise decreasing down to $u$.
\\ ${\ }$
In any almost complex manifold  $(X,J)$ each point has a fundamental neighborhood system of  $J$-pseudoconvex domains, and so 
the theorem above establishes local smooth approximation on $X$.
\\ ${\ }$
This result was proved in complex dimension 2 by the third author, who also showed that
the result would hold in general dimensions if a parallel result for continuous approximation were known.
This paper establishes the required step by solving the obstacle problem.

}}

\vskip .3in
\centerline{\bf Table of Contents}

\noindent
 \AA.     Introduction.   
 
\noindent
\BB.    The Obstacle Problem and Continuous Approximation in
for General Potential Theories.

\noindent
 \CC.    Strict Continuous Approximation of Plurisubharmonic Functions   on Almost
Complex Manifolds.

\noindent
 \DD.    Strict Smooth Approximation of Plurisubharmonic Functions   on Almost
Complex Manifolds.

 Appendix A.  Affine Jet-Equivalence.

 Appendix B.  $\Sigma_m$-Subharmonic Functions

\vfill\eject


\medskip
\noindent{\headfont \AA.\  Introduction.}

On any smooth almost  complex manifold $(X,J)$ there is a well-defined notion of $J$-plurisubharmonic
functions of class $C^2$, namely those $u\in C^2(X)$ which satisfy the condition $i\partial \dbar u \geq0$.
This notion  extends directly  to the space of distributions $\cd'(X)$ by requiring  the current
$i\partial \dbar u$ to be positive. It also extends to the space
$\USC(X)$ of upper semi-continuous functions $u:X\to [-\infty, \infty)$ in several ways -- using viscosity theory, 
or by requiring that  the restrictions   to J-holomorphic curves in $X$ be subharmonic. These different extensions have been shown to be, in a precise sense, equivalent (see [\Pa], [\AC]),
and the space of such functions is denoted by $\PSH(X,J)$.  

We say that a function $u\in C^2(X)$ is {\bf strictly} $J$-plurisubharmonic if $i\partial \dbar u >0$ at every point.
The manifold $X$ is then said to be {\bf  $J$-pseudoconvex} if it admits  a smooth (proper) exhaustion  function 
$\rho:X \to \bbr$ which is strictly $J$-plurisubharmonic. (See Remark \CC.7 for other equivalent definitions.)

The main point of this paper is to establish the following (in \S 4).

\noindent
{\bf THEOREM \DD.1. ($C^\infty$ Strict Approximation).} {\sl Suppose  $(X,J)$ is an almost complex
manifold which is $J$-pseudoconvex, and let $u \in \PSH(X,J)$ be a $J$-plurisubharmonic function.
Then there exists a decreasing sequence $\{u_j\} \ss C^\infty(X)$ of smooth strictly 
$J$-plurisubharmonic functions such that $u_j(x)\downarrow u(x)$ at each $x\in X$.}

Now  on  any almost complex manifold   $X$  every point  $x$
 has a fundamental neighborhood system of   $J$-pseudoconvex
domains   -- namely, small balls about $x$ in appropriate local coordinates.
Consequently, as a special case of Theorem  \DD.1 we have local  $C^\infty$ strict approximation
on $X$  (see Corollary \DD.2).

By this local regularization result a current $i\partial\bar{\partial} u\wedge i\partial\bar{\partial} v$ defined in
[\Pl] is a positive current for plurisubharmonic $u,v$ in the
Sobolev class $W^{1,2}_{loc}$, in particular for bounded
plurisubharmonic $u,v$ (see Proposition 4.2 and
Proposition 5.2 there and compare with Corollary 2 in [\Pll]).
For an application of our global regularization result see Corollary \DD.3, which concerns hulls of sets.

We note that in the case of plurisubharmonic functions on domains in $\bbc^n$,   smoothing as in Theorem  \DD.1
is  possible  on all pseudoconvex,  Reinhardt,  and tube domains (see [\FW]), but there are smooth domains where not all plurisubharmonic functions are a limits of a decreasing sequence of smooth plurisubharmonic functions (see [\Fo]).

Theorem  \DD.1 was proved in complex dimension 2 by the third  author (in [\Pll]), who
pointed out that his work would establish the result  in general dimensions provided one could prove
a certain parallel {\sl continuous} approximation theorem.  The required continuous approximation result
can be deduced from work of the first two  authors on the obstacle problem -- more precisely
the Dirichlet problem with an obstacle function.

The discussion of this obstacle problem in [\DDR] and [\Survey] and its exact implementation in 
the context of almost complex analysis is somewhat scattered, and so, for
clarity, we give a coherent exposition of the needed results
in the first two sections of this note. Nevertheless, this note draws heavily on the work
in  [\DDR], [\AC], [\Survey], [\Pl] and [\Pll].

It is interesting to note that the work in [\Pl] and [\Pll] also involves solving the Dirichlet
problem for the (almost) complex Monge-Amp\`ere operator. In this case, however,   the
solutions are taken in the smooth category using results in  [\Plis], where the
techniques are quite different from the viscosity methods employed in 
[\DDR], [\AC],  [\Survey]. The idea of using the   Monge-Amp\`ere equation
to approximate $J$-plurisubharmonic functions is probably due to J.-P. Rosay.

\noindent
{\bf Remark.}  The main proof in this paper consists of combining a Richberg-type
theorem (cf. [\Pl, Thm. 3.1], [\GPSH, Thm. 9.10])
 with the continuous approximation theorem
 which follows from solving the obstacle problem.  The method applies  generally to give smooth
approximation of $F$-subharmonic functions whenever these two components can be established.
An example  is given in Appendix B where smooth approximation is  established for
subsolutions of the complex Hessian equations on a K\"ahler  manifold.


\vfill\eject

\medskip
\noindent{\headfont \BB.\  The Obstacle Problem and Continuous Approximation for General Potential Theories.}

 We refer the reader to [\DDR] or [\Survey] for the concepts and terminology employed in this section. 
 
 Let $J^2(X)  \to X$ be the bundle of 2-jets of real-valued functions
 on a manifold $X$.  There is a natural splitting 
 $J^2(X) = \bbr \times J^2_{\rm red}(X)$ where the first factor corresponds
 to the value of the function.  
 
 Consider a subequation of the form $F = \bbr\times F_0$ with $F_0
 \ss J^2_{\rm red}(X)$. For a domain $\O \ss\ss X$,  let  $F(\ob)$ denote the set of 
 $u\in \USC(\ob)$ such that $u\bigr|_\O$ is $F$-subharmonic (i.e.,  $u\bigr|_\O$ is a  viscosity $F$-subsolution,
 cf. [\CRA], [\CIL]).

\Theorem{ \BB.1. (The Obstacle Problem)}
{\sl
Suppose that:

(1) $F_0$ is locally affinely jet-equivalent to a constant coefficient
(reduced) subequation $\bbf_0$,

(2)  $F_0$ has a monotonicity cone $M_0$ and $X$ carries a $C^2$ 
strictly $M$-subharmonic function $\psi$ where $M=\bbr\times M_0$,

(3) $g\in C(X)$, and

(4) $\O\ss\ss X$ is a domain with smooth boundary $\bo$ which is both
$F$- and $\ft$-strictly convex.

\noindent
Then the  function
$$
h(x) \ \equiv\ \sup_{u\in \cf[g]} u(x),
\eqno{(\BB.1)}
$$
where $\cf[g]\ \equiv\ 
\{u(x) : u\in F(\ob) \ {\rm and}\ u\leq g \ {\rm on}\ \ob\}$,  satisfies:

(i)\ \ \ $h\in C(\ob) \cap F(\ob)$,

(ii)\ \ $h\leq g$\ \ on $\ob$

(iii)\ $h\bigr|_{\bo} \ =\ g\bigr|_{\bo}$


\noindent
Furthermore,

(v) $h$ is the Perron function, and $\cf[g]$
is the Perron family, for the Dirichlet problem for the subequation 
$$
F^g \ \equiv (\bbr_- + g)\times F_0 \quad {\rm on}\ \O
$$
with boundary function $\vf\equiv g\bigr|_{\bo}$.

(vi)  Comparison holds for $F^g$ on $X$.
}

\vfill\eject

\noindent
{\bf COROLLARY \BB.2. (Continuous Strict Approximation).}
{\sl
Suppose $u\in F(\ob)$. 

 (a)  Then there exists a sequence  of functions $u_j \in C(\ob)\cap F(\ob)$
decreasing down to $u$ on $\ob$.  In fact, if $\{g_j\}\ss C(\ob)$ is any sequence
of continuous  functions decreasing down to $u$, the  $\{u_j\} \ss C(\ob)\cap F(\ob)$ can be chosen 
so that
$$
u\ \leq\ u_j\ \leq\ g_j \qquad \forall\, j.
\eqno({\BB.2)}
$$

(b) Moreover, given  $\e_j\downarrow 0$, the sequence $\{u_j + \e_j \psi\}$ 
also decreases down to $u$ on $\ob$, and on each compact subset
of $\O$, the functions $\{u_j + \e_j \psi\}$  are $c$-strict for some $c>0$.
}

See  \BB.3 below for a definition and discussion of $c$-strictness.

\noindent
{\bf Proof of Corollary \BB.2.}  Pick $g_j \in C(\ob)$ with $g_j \downarrow u$.  Let $u_j$ denote 
the solution of the obstacle problem for $g_j$.  Then  $u_j \in C(\ob)\cap F(\ob)$
and $u_j\leq g_j$.  Since $u$ is in the Perron family $\cf[g_j]$, we have (\BB.2).
 This proves Part (a).  Part (b) follows from (a) and hypothesis (2).\qed

\medskip
\noindent
{\bf Proof of Theorem \BB.1. } 
The following is proved in [\DDR] but not stated explicitly as a theorem.
It is however stated explicitly as Theorem 8.1.2 in [\Survey] and the 
proof is given there based on results in [\DDR]

\Theorem{8.1.2 in [\Survey]}  {\sl Suppose $F$ is a subequation on a manifold $X$ which is
locally affinely jet-equivalent to a constant coefficient subequation.
Suppose there exists a $C^2$ strictly $M$-subharmonic function on $X$ where
$M$ is a monotonicity cone for $F$.
Then for every domain $\O\ss\ss X$
whose boundary is strictly $F$- and 
$\ft$-convex, both existence and uniqueness hold for the Dirichlet problem.
That is, for every $\vf \in C(\bo)$  there exists a unique $F$-harmonic function
$u\in C(\ob)$ with $u\bigr|_{\bo} =\vf$.}

The adaptation to the general Obstacle Problem is given in Section 8.6 of
[\Survey]. What follows is a more detailed version of that argument.

By assumption we  know that $F=\bbr\times F_0$ is affinely jet equivalent 
to the constant coefficient equation $\bbr\times \bbf_0 \ss  \bbr \times \bbj^2_{\rm red}$,
with a jet equivalence which is the identity on the first factor.
Hence the subequation 
$$
F^g \  \equiv \ \{r\leq g(x)\}\times F_0
$$
is locally affinely jet equivalent to the subequation
$$
\bbf^g \ \equiv\ \{r\leq g(x)\}\times \bbf_0
$$
We now consider the affine jet equivalence
$$
\Phi : \bbr\times \bbj^2_{\rm red} \ \arr\  \bbr\times \bbj^2_{\rm red}
$$
given by
$$
\Phi(r, J) \ \equiv\ (r-g(x), J).
$$
  Applying this gives the local equivalence
  $$
  \Phi : \bbf^g \ \arr\ \{r\leq0\}\times \bbf_0 \ \equiv \ \bbr_- \times \bbf_0,
  $$
and so composing this with the first equivalence shows that $F^g$ is locally
affinely jet-equivalent to the constant coefficient subequation $\bbr_- \times \bbf_0$.

Now observe that if $M_0$ is a monotonicity cone for $F_0$,
then $M_- \equiv  \bbr_-\times M_0$ is a monotonicity cone for $F^g$.

Note also that if $\psi$ is strictly $M$-subharmonic function, then so is $\psi - c$
for any constant  $c\leq 0$ because $M$ satisfies the  basic negativity condition (N). 
Given a domain $\O\ss\ss X$, we may therefore
assume that $\psi <0$ on a neighborhood of $\ob$.  In this case, $\psi$ is also
$M_-$-strictly subharmonic on $\ob$.

Since $F^g$ is locally jet-equivalent\begin{footnote} 
{See Appendix A for a discussion of jet-equivalence.}
\end{footnote}
to a constant coefficient subequation,
local weak comparison holds for $F^g$.  This is Theorem 10.1 in [\DDR] 
and follows from the   Theorem on Sums.
Local weak comparison implies weak comparison (Theorem 8.3 in [\DDR]).
Now using Theorems 9.5 and 9.2 we have that comparison holds for $F^g$ on $X$.

The Dirichlet Problem for $F^g$-harmonics would now be solvable for 
arbitrarily prescribed boundary data $\vf\in C(\bo)$,  (by either Theorem 12.4 in [\DDR] 
or  Theorem 8.1.2 above)  if one   could  prove
that the boundary is strictly $F^g$ and $\wt {F^g}$ convex.  

However, this is not true in general, and in fact 
existence fails for a boundary function $\vf\in C(\bo)$ unless 
 $\vf \leq g\bigr|_{\bo}$.  Nevertheless, {\sl   if $\bo$ is both $F$ and $\wt F$ strictly convex, then existence holds for each boundary function $\vf \leq g\bigr|_{\bo}$. }
Section 8.6 in [\Survey] provides a  proof of this.

Here we give a proof 
but with attention restricted  to the case
 at hand where $\vf = g\bigr|_{\bo}$. 
The Perron family for $F^g$ with this boundary
 data consists of those  functions $u\in\USC(\ob)$ which are 
 $F$-subharmonic on $\O$ and satisfy the  additional constraint that
$u\leq g$ on $\O$. The dual subequation to $F^g$ is 
$\wt {F^g} = [(\bbr_--g)\times J^2_{\rm red}(X) ] \cup \ft$.  
Since $\wt {F^g} \ss  \ft$, the $\bo$ is strictly
$\wt  {F^g} $-convex if  it is strictly $\ft$-convex.  However,  $\bo$ can never be strictly
$F^g$-convex, as defined in Definition 11.10 of [\DDR],  because  $(\oa{{F_\la}})_x = \emptyset$ for $\la> g(x)$),

Nevertheless, the only place that this hypothesis is used 
in proving  Theorem 8.1.2 for $H$ is in the barrier construction which appears in the proof of 
 Proposition $F$ in [\DDR].
With $\vf(x_0) = g(x_0)$,  the barrier $\b(x)$ as defined in (12.1) in [\DDR]
 is not only $F$-strict near $x_0$ but also automatically $F^g$-strict since $\b < g$
 in a neighborhood of $x_0$.  \qed

\noindent
{\bf Definition \BB.3. (Strictness).}  Let $F\ss J^2(X)$ be a subequation.     A function $u\in F(\O)$ is {\bf strictly $F$-subharmonic}  (or simply {\bf  strict}) 
 if for any $\vf \in  C_0^\infty(\O)$, there exists $\e>0$ such that $u+\e\vf \in F(\O)$.
 
Note that a  $C^2$-function $u\in F(\O)$ is strict iff $J^2_x u\in \Int F$ $\forall\,x\in\O$.

In [\DDR] there is the following  related concept of $c$-strictness for $c>0$.
 Equip $J^2(X)$ with a bundle metric  (induced, say, 
from a riemannian metric on $X$), and for $x\in X$, define $F_x^c\equiv \{J \in F_x : \dist_x(J, \sim F )\geq c\}$
where $\dist_x$ denotes the distance in the fibre.  A function $u\in F(\O)$ is said to be {\bf $c$-strict}
 on a compact set $K\ss  \O$ if $u$ is $F^c$-subharmonic on a neighborhood of $K$. The constant
 $c$ depends on the choice of bundle metric, but the condition of being $c$-strict on $K$ for some
 $c>0$ does not.  Strictness, as defined above, is equivalent to being locally $c$-strict on $\O$.
 (This is proved, though not explicitly stated, in \S 7 of [\DDR].)
 
 \noindent
 {\bf Remark \BB.4.}  The main conclusion of Theorem \BB.1 above can be stated in more appealing 
 and succinct terms.  Let us call  the function $h$, defined in (\BB.1), the {\bf largest $F$-subharmonic
 minorant of $g$}.  Then we have the following abbreviated version of Theorem \BB.1 and Corollary \BB.2.  

\noindent
{\bf THEOREM  \BB.5.} {\sl  Suppose $X,F= \bbr\times F_0$ and $\O$ are as in Theorem \BB.1.
Then given $g\in C(\ob)$, the largest $F$-subharmonic  minorant of $g$ on $\ob$ is continuous
and equals $g$ on the boundary of $\O$.
 \\ ${\ }$
Moreover, given $u\in F(\ob)$ there exists a sequence $\{u_j\} \ss C(\ob) \cap F(\ob)$ decreasing 
down to $u$ (with each $u_j$ strict).
}

\vskip .3in


\noindent
{\headfont \CC. Strict Continuous Approximation of Plurisubharmonic Functions  
 on Almost
Complex Manifolds}

Let $(X,J)$ be an almost complex manifold, and let $F(J) \ss J^2_{\rm red}(X)$
be the subequation defining the upper semi-continuous $J$-plurisubharmonic functions
on $X$.  (It is shown in [\AC] that all the different  basic definitions of these functions
are, in a precise sense, equivalent).\begin{footnote}{It is also shown at the end of section 7 in [\AC]
that the various notions of $F(J)$-harmonic (including the notion of being  maximal and continuous) 
are equivalent.}\end{footnote}

Proposition 4.5 in the paper  [\AC] proves that the subequation
$F(J)$ is locally jet equivalent to  a constant coefficient reduced subequation  (in fact to the standard
 subequation $F(J_0)\cong\{i\partial\dbar u\geq0\}$ determined by a standard parallel $J_0$).

Furthermore,  $F(J)$ is a convex cone subequation and in particular it
satisfies $F(J) + F(J)\ss F(J)$.  Therefore, $F(J)$ is a monotonicity cone
for itself.  A $C^2$-function $\psi$ is  strictly $J$-plurisubharmonic  
(i.e.,  strictly $F(J)$-subharmonic)  if
$i\partial\dbar \psi >0$ on $X$.  

\noindent
{\bf Definition \CC.1.}  A domain $\O\ss\ss X$ is called {\bf strictly $J$-pseudoconvex} if it has a global $C^2$
defining function $\psi$ which is strictly $J$-plurisubharmonic  
 on a neighborhood of $\ob$. Let $\wt F(J)$ denote the dual subequation.
One  checks that 
$$
F(J) + F(J)\ss F(J) \ \ \Rightarrow\ \  \wt F(J) + F(J)\ss \wt F(J)
\ \ \Rightarrow\ \ F(J)\ss \wt F(J),
$$
so if $\bo$ is strictly $F(J)$-convex, it 
is automatically strictly $\wt F(J)$-convex.

Thus,  as a special case of Theorem \BB.5  we have the following.

\noindent
{\bf  THEOREM \CC.2.} {\sl Let $\O\ss\ss X$ be a strictly $J$-pseudoconvex domain
in an almost complex manifold $(X,J)$. Let $g \in C(\ob)$. Then the largest
$J$-plurisubharmonic minorant of $g$ is continuous.
\\ ${\ } $
Moreover, given $u\in \PSH(\ob)$ there exists a sequence $\{u_j\}\ss C(\ob)\cap \PSH(\ob)$ decreasing down
to $u$ (with each $u_j$ strict).
}


We now address the global question of continuous approximation of 
$J$-plurisubharmonic functions on $X$.

\noindent
{\bf Definition \CC.3.}  An almost  complex manifold  $(X,J)$ is  {\bf $J$-pseudoconvex} if it has a global $C^2$ strictly $J$-plurisubharmonic exhaustion function. 
 (See Remark \CC.7 below for equivalent definitions.)

It is standard that a strictly $J$-pseudoconvex domain $\O$ is itself $J$-pseudo-convex.

\noindent
{\bf  THEOREM \CC.4.} {\sl Suppose  $X$ is a  $J$-pseudoconvex manifold. Then for each $u\in{\rm PSH}(X)$
there  exists a sequence of continuous strictly $J$-plurisubharmonic functions $u_j \in C(X)$ 
decreasing down to $u$ on $X$.}

\noindent
{\bf Proof.}
We  shall adapt a part of the proof of the Theorem 1 from [\Pll].  
Take a decreasing sequence of continuous functions $\{g_k\}$ converging down  to $u$.
We begin with a result in smooth topology.

\noindent
{\bf Claim \CC.5.} Let $h$ be an arbitrary continuous function on $X$, and
suppose that  $\rho:X \to \bbr$  is a  $C^2$ (proper) exhaustion function.  
Then there exists a convex function $\chi\in C^\infty(\bbr)$ with $\chi'  \geq 1$
so that 
$$
\chi ( \rho(x) ) \ \geq \ h(x) \qquad\text{for all $x\in X$}.
$$

\noindent
{\bf Proof.}
Set $\psi(t) \equiv \sup\{h(x) : \rho(x) \leq t\}$ and note that 
$$
\chi(\rho(x)) \ \geq\ h(x) \ \  \forall\, x\in X 
\qquad\iff\qquad \chi(t) \ \geq\ \psi(t) \ \  \forall\, t\in {\rm range}(\rho). 
$$
This reduces the claim to a one-variable claim.
To establish this,  assume that  range$(\rho) = [0,\infty)$ and replace 
$\psi$ by a smooth function which is larger.  Then choose $\chi \in C^\infty([0,\infty))$
to have $\chi(0) = \psi(0)$,
$\chi'(0) \geq \max\{\psi'(0), 1\}$ and  $\chi'' \geq \max\{\psi'', 0\}$. \qed

 Now let $\rho \in C^\infty(X)$ be a  strictly $J$-plurisubharmonic exhaustion function.
For any smooth convex, increasing function $\chi\in C^\infty(\bbr)$, with $\chi' \geq 1$, the composition
$\chi\circ \rho$ is also  a smooth  strictly $J$-plurisubharmonic exhaustion.
Thus, by Claim \CC.5, with $h$ taken to be $g_1$ plus any exhaustion function for $X$,
 we can assume  $\rho$ is chosen so that 
$$
\lim_{z\rightarrow\infty}(\rho(z)-g_1(z))=+\infty
\eqno{(\CC.1)}
$$
where $\lim_{z\rightarrow\infty}$ denotes the limit in the one-point compactification of $X$.

By (\CC.1) the  sets $U_k \equiv \{\rho > g_1+k\}$ provide a fundamental neighborhood system for the point 
at infinity.  Since $\rho$ is an exhaustion, we have that  $\{\rho - k \geq t\} \ss U_k$ if $t$ is sufficiently large.
By Sard's Theorem we may choose such  $t$ to be a regular value $t_k$ of $\rho-k$.
Then $\O_k \equiv \{\rho-k<t_k\}$ is a strictly $J$-pseudoconvex domain, and 
$$
\rho - k\ >\  g_1 \ (\geq \ g_k) \quad \text{on a neighborhood of }\  \sim \O_k.
\eqno{(\CC.2)}
$$
Hence, 
$$
\wt{g}_k \ \equdef\ \max\{g_k, \rho-k\} \   = \  \rho-k
\quad \text{on a neighborhood of }\  \sim \O_k.
\eqno{(\CC.3)}
$$

Now let  $u_k$ be the  largest $J$-psh minorant
 of $\tilde{g}_k$ on $\O_k$, and note that  $u_k$ is continuous by  Theorem \CC.2.  
By (\CC.3) we have  $\tilde{g}_k  = \rho - k$ on a  neighborhood of  $\sim \O_k$.
Since $ \rho - k$  is $J$-psh, 
and $u_k$ is the largest $J$-psh minorant of $\tilde{g}_k$, we have 
$u_k=\rho-k$ on  a  neighborhood of  $\sim \O_k$.
 Thus  we can extend $u_k$  as a $J$-psh 
function to  all of $X$ by setting $u_k  =  \rho-k$ on $\sim \O_k$.

Note that since $\wt {g}_{k} \equiv   \max\{g_{k}, \rho-k\}$, 
$g_{k+1} \leq g_k$, and $g_k \downarrow u$, one has 
$$
\wt {g}_{k+1} \ \leq\  \wt{g}_k   \and   \wt{g}_k \ \downarrow\ u.
\eqno{(\CC.4)}
$$
By definition
$$
u_{k} \leq \wt{g}_{k} 
\and 
u_{k} = \wt{g}_{k} \ \ {\rm on} \  \sim \O_k.
\eqno{(\CC.5)}
$$
Now since $u_{k+1} \leq \wt{g}_{k+1}$, and since $u_k$ is the
largest $J$-psh minorant of $\wt{g}_k$ on $\ob_k$, we have by (\CC.4) that 
$u_{k+1} \leq u_k$ on $\ob_k$.  On the complement $\sim \O_k$, we
have $u_k = \wt{g}_k$ and so  $u_{k+1} \leq u_k$ again  by (\CC.4) and (\CC.5).
Hence,
$$
u_{k+1} \ \leq\ u_k \qquad{\rm on} \ X.
\eqno{(\CC.6)}
$$
Since $u\leq \wt{g}_k$ is $J$-psh and $u_k$ is the largest such minorant 
on $\ob_k$, we have that
$u \leq u_k$ on  $\ob_k$.
On the complement $\sim \O_k$, we
have $u_k = \wt{g}_k$ and so  $u \leq u_k$ there as well.
Hence,
$$
u \ \leq\ u_k  \and  u_k\ \downarrow\ u \qquad{\rm on} \ X.
$$
In other words $\{u_k\}$  is a  decreasing sequence 
of continuous $J$-psh functions  decreasing down  to $u$ on $X$, and we can replace  $u_k$  
with  $u_k+\frac{1}{k}\rho$   to make $u_k$ strict. \qed


\Remark{\CC.7. (Equivalent Definitions of $J$-Pseudoconvexity)} 
In defining $J$-pseudoconvexity it is  enough to assume the existence of a {\sl continuous}
 strictly $J$-plurisubharmonic exhaustion function $\rho:X\to \bbr$. This follows from the 
 extension of Richberg's Theorem  to  almost complex manifolds  (Theorem 3.1 in [\Pl]).
 Such manifolds are called   {\sl almost Stein manifolds} in  [\DS].

$J$-Pseudoconvex manifolds $(X,J)$ can also be characterized in terms of the hulls of compact sets
(see (\DD.1) below) by  requiring that:
\\ ${\ }$ \qquad
(i) \  There exists some $u\in \PSH^\infty(X,J)$ which is strict, and 
\\ ${\ }$ \qquad
(ii)  For every compact $K\ss X$, the hull $ {\wh K}_{C^\infty}$ is compact.
 \\
 By Theorem 3.1 in [\Pl] we have that the hulls  $\wh{K}_{C^0} = \wh{K}_{C^\infty}$ agree
 (see Corollary \DD.3 below).
Therefore, $J$-Pseudoconvex manifolds can also be characterized by the requiring:
\\ ${\ }$ \qquad
(i) \  There exists some $u\in \PSH^0(X,J)$ which is strict, and 
\\ ${\ }$ \qquad
(ii)  For every compact $K\ss X$, the hull $ {\wh K}_{C^0}$ is compact.
 \\
For the proof one applies  standard arguments  (cf. [\GPSH, \S4] or [\PUP, Prop. 9.3])
to show that (i) and (ii) imply the existence of a strict PSH-exhaustion (in either case).

\vskip.3in

\noindent
{\headfont \DD. Strict Smooth Approximation of Plurisubharmonic Functions  
 on Almost   Complex Manifolds}

\noindent
{\bf THEOREM \DD.1. ($C^\infty$ Strict Approximation).}  {\sl Suppose  $(X,J)$ is an almost complex
manifold which is $J$-pseudoconvex, and let $u \in \PSH(X,J)$ be a $J$-plurisubharmonic function.
Then there exists a decreasing sequence $\{u_j\} \ss C^\infty(X)$ of smooth strictly 
$J$-plurisubhrmonic functions such that $u_j(x)\downarrow u(x)$ at each $x\in X$.}

\noindent
{\bf Proof.} Apply  Theorem 3.1 in [\Pl] and Theorem \CC.4 above.\qed

This  generalizes Theorem 1 in [\Pll] to arbitrary dimensions.

\noindent
{\bf COROLLARY \DD.2. (Local  $C^\infty$ Strict Approximation).} {\sl 
Let $(X,J)$ be an arbitrary (smooth) almost complex manifold.
Then every point $x \in X$ has a fundamental system of neighborhoods
 $U$ with the property that for every $u\in\PSH(U,J)$ there is a decreasing sequence 
$\{u_j\}\ss  C^\infty(U)$ of strictly $J$-plurisubharmonic functions such that $u_j\downarrow u$.}

\noindent
{\bf Proof.}  Fix local coordinates in $\bbc^n$ for $X$ near $x$ so that $J$ is $C^1$-close to
the standard $J_0$ at the origin.  Then $\chi(z) = |z|^2$ is strictly $J$-psh on the ball 
$B_{\e}(0)  = \{|z| < \e\}$ for all $\e>0$ sufficiently small.
It is standard that any domain which admits a $C^2$ strictly $J$-plurisubharmonic defining 
function, is $J$-pseudoconvex.
\qed

One can also give a more direct proof of Corollary \DD.2 based on Theorem \CC.2 above
and Theorem 3.1 in [\Pl].

Another  immediate consequence of the global approximation Theorem \DD.1 is that all the
various possible definitions of the hull of a set actually agree.
Given a compact set  $K \ss X$ we define its {\bf $J$-plurisubharmonic hull} to be the set
 $$
 \wh K \ \equiv\ \left\{x\in X :  u(x) \leq \sup_K u \ \ \forall\, u\in \PSH(X,J)\right\}.
\eqno{(\DD.1)}
$$
 One could also define ${\wh K}_{C^0}$  and ${\wh K}_{C^\infty}$ by replacing 
 $\PSH(X,J)$ in (\CC.4) with  $\PSH^0(X,J) \equiv \PSH(X,J)\cap C(X)$ and  
 $\PSH^\infty(X,J) \equiv\ \PSH(X,J)\cap C^\infty(X)$
 respectively.

\noindent
{\bf Corollary \DD.3.} {\sl Suppose $(X,J)$ is $J$-pseudoconvex.
Then for any  compact $K\ss X$,   one has} 
${\wh K} \  =\  {\wh K}_{C^0} \ =\  {\wh K}_{C^\infty}$.

\noindent
{\bf Proof.}  Clearly ${\wh K} \ss  {\wh K}_{C^0} \ss {\wh K}_{C^\infty}$, so it suffices 
to show that $ {\wh K}_{C^\infty} \ss  {\wh K}$.
Suppose that $x\notin  {\wh K}$.  Then there exists $u\in \PSH(X,J)$ with $u\leq0$ on $K$ and  $u(x)=1$.
Replace $u$ with $\max\{u, 0\}$.  Let $\{u_j\}$ be the sequence given in Theorem \DD.1.  Then $u_j$ converges uniformly to 0 on the compact set $K$ 
and $u_j(x) \geq1$ for all $j$.  Hence, $x\notin {\wh K}_{C^\infty}$. \qed

\vskip.3in

\noindent
{\headfont Appendix A. Affine Jet-Equivalence.}
A local affine jet-equivalence is a local isomorphism of the 2-jet bundle $\bbj(\rn) = \bbr\times
\rn\times\Symn$
which is of the form:
$$
r' = r+r_0(x), \quad p'=k(x)p + p_0(x),\quad A'=    h(x) A h(x)^t + L_x(p) +A_0(x)
$$
where
$$
\begin{aligned}
&r_0(x) \ \text{takes values in $\bbr$}, \\
&p_0(x) \ \text{takes values in $\rn$}, \\
&A_0(x) \ \text{takes values in $\Symn$}, \\
&\text{(i.e., $J_0(x) \equiv (r_0(x), p_0(x), A_0(x))$ is a section of $\bbj(\rn)$)}
\end{aligned}
$$
and 
$$
\begin{aligned}
&\text{ $k(x)$ and $h(x)$ take values in ${\rm GL}_n(\bbr)$, while} \\
&\text{ $L_x$ takes values in $\Hom(\rn, \Symn)$}\\
\end{aligned}
$$

The regularity conditions on the jet-equivalence required in the proof of 
Theorem 10.1 in [\DDR] are:
\smallskip

(1) \ \ $k,h$ and $L$ are Lipschitz continuous, and

\smallskip

(2) \ \ $J_0$ is continuous.
\medskip

For the second jet equivalence in our application to the 
Obstacle Problem, $g\equiv h \equiv Id$ and $J_0(x) = (r_0(x),0,0)$,
so our obstacle function $g(x) = -r_0(x)$ need only be continuous.


\vskip.3in


\noindent
{\headfont Appendix B.   $\Sigma_m$-Subharmonic Functions.}

As noted in Remark \AA.3, for any subequation $F$, smooth approximation for 
$F$-subharmonic functions can be proved
whenever continuous approximation and a Richberg-type theorem can be established for $F$.
In this appendix we give just such a result for the complex hessian subequations on a K\"ahler manifold.

 Let $X$ be a complex manifold of   dimension $n$ with a fixed K\"{a}hler form $\omega$. We say that a function $u\in\mathcal{C}^2(\Omega)$ is $\Sigma_m$-subharmonic on a domain 
$\Omega\ss\ss X$ if $(dd^cu)^k\wedge\omega^{n-k}\geq0$ for $k=1,\ldots,m$. 
We say that a locally integrable function $$u:\Omega\rightarrow[-\infty,+\infty)$$ is $\Sigma_m$-subharmonic 
($u\in {\Sigma}_m(\Omega)$)
 if $u$ is upper semicontinuous and $$dd^cu\wedge dd^cu_1\wedge\ldots\wedge dd^cu_{m-1}\wedge\omega^{n-m}\geq0,$$ for any  $\mathcal{C}^2$  $\Sigma_m$-subharmonic functions $u_1,\ldots u_{m-1}$ (they are defined in [\BLO] for $\omega=\omega_{st}=dd^c(|z|^2)$ in $\mathbb{C}^n$ and in [\DK]  and [\Lu] for general K\"{a}hler form).  
This is just the subequation $F\equiv \Sigma_m$ defined on $X$ by the condition that the first $m$ elementary
symmetric functions of the complex hessian satisfy $\sigma_\ell(\Hess_\bbc  u) \geq 0$ for $\ell=1,...,m$
(compare Example 18.1 in [\DDR] and Lemma 7 in  [\Plll]).

A Richberg-type theorem for ${\Sigma}_m$  was proved in [\Plll] (Theorem 2).
 Lu and Nguyen proved in [\LN] that on compact K\"{a}hler manifolds any quasi-$\Sigma_m$-subharmonic function
can be approximated from above by smooth quasi-$\Sigma_m$-subharmonic functions
 (a function $u$ is quasi-$\Sigma_m$-subharmonic if  the  function $u+\rho$ is $\Sigma_m$-subharmonic where $\rho$ is local potential for $\omega$). Actually their global result implies that locally it  is possible to regularize $\Sigma_m$-subharmonic functions.  However, in the same way as in Theorem 4.1, we can prove a slightly stronger result.

\Theorem{B.1} {\sl
Let $ X$ be a  $\Sigma_m$-pseudoconvex K\"ahler manifold. Let $u$ be a $\Sigma_m$-subharmonic function
on $X$. Then there exists a decreasing sequence $u_j \in  \mathcal{C}^\infty(X)$ of 
$\Sigma_m$-subharmonic functions such that $u_j\downarrow u$.
}

By {\bf $\Sigma_m$-pseudoconvex} we mean that  $X$ has a global $\mathcal{C}^2$ strictly
$\Sigma_m$-subharmonic exhaustion function. In particular Stein manifolds are $\Sigma_m$-pseudoconvex.




\def\item{}
\vfill\eject

\centerline{\headfont References}

\noindent
\item{[\BLO]}  
 Z. B\l ocki, \textit{Weak solutions to the complex Hessian
   equation}, Annales de l'Institut Fourier 55 (2005), 1735-1756.
 \smallskip

\noindent
\item{[\CRA]}   M. G. Crandall,  {\sl  Viscosity solutions: a primer},  
pp. 1-43 in ``Viscosity Solutions and Applications''  Ed.'s Dolcetta and Lions, 
SLNM {\bf 1660}, Springer Press, New York, 1997.

 \smallskip

\noindent
\item{[\CIL]}   M. G. Crandall, H. Ishii and P. L. Lions {\sl
User's guide to viscosity solutions of second order partial differential equations},  
Bull. Amer. Math. Soc. (N. S.) {\bf 27} (1992), 1-67.

\smallskip

\noindent
\item{[\DS]} 
K. Diederich and A. Sukhov,
{\sl Plurisubharmonic exhaustion functions and almost complex Stein structures},
Michigan Math. J. {\bf 56} (2008), no. 2, 331-355.

 \smallskip

\noindent
\item{[\DK]}
S. Dinew and  S. Ko\l odziej,   {\sl A priori estimates for complex Hessian equations}, Anal. PDE 7 (2014), no. 1, 227-244.

 \smallskip

\noindent
\item{[\Fo]}  J. E. Forn\ae ss,  {\sl
Plurisubharmonic functions on smooth domains},
Math. Scand. {\bf 53} (1983), no. 1, 33-38.

 \smallskip

\noindent
\item{[\FW]}  J. E. Forn\ae ss and J. Wiegerinck, {\sl Approximation of plurisubharmonic functions}, 
Ark. Mat. {\bf 27} (1989), no. 2, 257-272.

\smallskip

 \noindent
\item{[\DDD]}   F. R. Harvey and H. B. Lawson, Jr., {\sl  Dirichlet duality and the non-linear Dirichlet problem},    Comm. on Pure and Applied Math. {\bf 62} (2009), 396-443. ArXiv:math.0710.3991.

\smallskip

 \noindent
\item {[\PUP]}   \ \----------,  {\sl  Plurisubharmonicity in a general geometric context},
  Geometry and Analysis {\bf 1} (2010), 363-401. ArXiv:0804.1316

\smallskip

 \noindent
\item{[\DDR]}  \ \----------, {\sl Dirichlet Duality and the Nonlinear Dirichlet Problem on Riemannian Manifolds},  J. Diff. Geom. {\bf 88} (2011), 395-482.   ArXiv:0912.5220.

  \smallskip

 \noindent
\item {[\GPSH]} \ \----------,  {\sl  Geometric plurisubharmonicity and convexity - an introduction},
  Advances in Math.  {\bf 230} (2012), 2428-2456.   ArXiv:1111.3875.

\smallskip

 \noindent
\item{[\AC]}  \ \----------,  {\sl  Potential theory on almost complex manifolds},  {Ann. Inst. Fourier}  (to appear).
ArXiv: 1107.2584.

\smallskip

 \noindent
\item {[\Survey]}  \ \----------,   {\sl  Existence, uniqueness and removable singularities
for nonlinear partial differential equations in geometry},\ 
 pp. 102-156 in ``Surveys in Differential Geometry 2013'', vol. 18,  
H.-D. Cao and S.-T. Yau eds., International Press, Somerville, MA, 2013.
ArXiv:1303.1117.

\smallskip

 \noindent
\item {[\Lu]} 
 H. C. Lu,  {\sl Solutions to degenerate complex Hessian equations}, J. Math. Pures Appl. (9) 100 (2013), no. 6, 785-805.

\smallskip

 \noindent
\item {[\LN]} 
H. C. Lu and V.-D. Nguyen, 
{\sl Degenerate complex Hessian equations on compact K\"{a}hler manifolds},  arXiv:1402.5147.

\smallskip

 \noindent
\item{[\Pa]}  \ \  N. Pali, {\sl Fonctions plurisousharmoniques et courants positifs de type $(1, 1)$ sur une vari\'{e}t\'{e} presque  complexe},
Manuscripta Math. {\bf 118} (2005), no. 3, 311-337.

\smallskip

 \noindent
\item{[\Plis]}  \ \  S. Pli\'s,  {\sl The Monge-Amp\`ere equation on almost complex manifolds},  Math. Z.
{\bf 276} (2014), no. 3-4, 969-983.

 \noindent
\item{[\Pl]}   \ \----------,  {\sl  Monge-Amp\`ere operator on four dimensional  almost complex manifolds},   
ArXiv: 1305.3461.

 \smallskip

\noindent
\item {[\Pll]}  \ \----------,   {\sl On regularization of
$J$-plurisubharmonic functions}, C.R. Acad. Sci. Paris, Ser.I (2014),
http://dx.doi.org/10.1016/j.crma.2014.11.001. 

\smallskip
 
 \noindent
\item{[\Plll]}  \ \----------,   {\sl The smoothing of $m$-subharmonic functions},
 arXiv:1312.1906.

 \end{document}